\documentclass[11pt]{amsart}
\usepackage{latexsym,amssymb,amsfonts,amsmath,graphicx,fullpage,url}
\usepackage[vcentermath,enableskew]{youngtab}
\usepackage{color}
\usepackage{cite}
\usepackage[all]{xy}
\usepackage{verbatim}
\usepackage{tikz,bbm}
\usepackage{caption}
\usepackage{subcaption}
\usepackage{hyperref}
\newtheorem{theorem}{Theorem}[section]
\newtheorem{proposition}{Proposition}[section]
\newtheorem{lemma}[theorem]{Lemma}
\newtheorem{corollary}[theorem]{Corollary}
\theoremstyle{definition}
\newtheorem{definition}[theorem]{Definition}

\theoremstyle{remark}

\numberwithin{equation}{section}

\def\C{{\mathcal{C}}}
\def\S{{\mathcal{S}}}

\def\P{\mathcal{P}}

\def\xpi{\overline{X_\pi}}
\def\torb{\overline{T^{2n}\cdot p}}
\def\torbq{\overline{T^{2n}\cdot q}}

\def\mpy{\Phi(\mathbb{P}(Y_\pi))}
\def\L{\overline{L(\pi)}}
\newcommand{\sizeQ}{m}

\newcommand{\be}{\begin{equation}}
\newcommand{\ee}{\end{equation}}

\newcommand{\bd}{\begin{definition}}
\newcommand{\ed}{\end{definition}}
\newcommand{\bt}{\begin{theorem}}
\newcommand{\et}{\end{theorem}}
\newcommand{\bl}{\begin{lemma}}
\newcommand{\el}{\end{lemma}}
\newcommand{\bp}{\begin{proposition}}
\newcommand{\ep}{\end{proposition}}
\newcommand{\bc}{\begin{corollary}}
\newcommand{\ec}{\end{corollary}}

\newcommand{\lmp}[1]{\mathrel{\mathop{\longmapsto}^{\mathrm{#1}}}} 

\newtheorem*{theorem1*}{Theorem \ref{thm:main}}

\newcommand{\old}[1]{}

\font\co=lcircle10
\def\petit#1{{\scriptstyle #1}}

\def\jr{\smash{\raise2pt\hbox{\co \rlap{\rlap{\char'005} \char'007}}
               \raise6pt\hbox{\rlap{\vrule height6.5pt}}
               \raise2pt\hbox{\rlap{\hskip4pt \vrule height0.4pt depth0pt
                width7.7pt}}}}
\def\je{\smash{\raise2pt\hbox{\co \rlap{\rlap{\char'005}
                \phantom{\char'007}}}\raise6pt\hbox{\rlap{\vrule height6pt}}}}
\def\+{\smash{\lower2pt\hbox{\rlap{\vrule height14pt}}
                \raise2pt\hbox{\rlap{\hskip-3pt \vrule height.4pt depth0pt
                width14.7pt}}}}
\def\perm#1#2{\hbox{\rlap{$\petit {#1}_{\scriptscriptstyle #2}$}}%
                \phantom{\petit 1}}

\def\textcross{\ \smash{\lower4pt\hbox{\rlap{\hskip4.15pt\vrule height14pt}}
                \raise2.8pt\hbox{\rlap{\hskip-3pt \vrule height.4pt depth0pt
                width14.7pt}}}\hskip12.7pt}
\def\textelbow{\ \hskip.1pt\smash{\raise2.8pt%
                \hbox{\co \hskip 4.15pt\rlap{\rlap{\char'005} \char'007}
                \lower6.8pt\rlap{\vrule height3.5pt}
                \raise3.6pt\rlap{\vrule height3.5pt}}
                \raise2.8pt\hbox{%
                  \rlap{\hskip-7.15pt \vrule height.4pt depth0pt width3.5pt}%
                  \rlap{\hskip4.05pt \vrule height.4pt depth0pt width3.5pt}}}
                \hskip8.7pt}

\title{Toric matrix Schubert varieties and their polytopes}
\author{Laura Escobar}
\address{Laura Escobar,
Institut f\"ur Mathematik, TU Berlin, Str.  des 17.
Juni 136, 10623 Berlin, Germany  \newline le78@cornell.edu
}

\author{Karola M\'esz\'aros}
\address{Karola M\'esz\'aros, Department of Mathematics, Cornell University, Ithaca NY 14853  \newline karola@math.cornell.edu
}
\thanks{M\'esz\'aros was partially supported by a National Science Foundation Grant  (DMS 1501059).}
\rightmark

\begin{document}

\maketitle

\begin{abstract} Given a  matrix Schubert variety $ \xpi$,  it can be written as $\xpi=Y_\pi\times \mathbb{C}^q$ (where $q$ is maximal possible). We characterize when $Y_{\pi}$ is toric (with respect to a $(\mathbb{C}^*)^{2n-1}$-action) and study the associated polytope $\Phi(\mathbb{P}(Y_\pi))$ of its projectivization. We construct regular triangulations of $\Phi(\mathbb{P}(Y_\pi))$ which we show are geometric realizations of a family of subword complexes. Subword complexes were introduced by Knutson and Miller in 2004, who also showed that they are  homeomorphic to balls or spheres and raised the question of their polytopal realizations. 
\end{abstract}  

\section{Introduction}
In this paper we study the  geometry of matrix Schubert varieties and use it to give geometric realizations of a family of subword complexes. Matrix Schubert varieties were introduced in \cite{MR1154177} to study the degeneraci loci of flagged vector bundles. In \cite{annals} Knutson and Miller showed that Schubert polynomials are multidegrees of matrix Schubert varieties.  
On the other hand, in \cite{subword} Knutson and Miller introduced subword complexes  to illustrate the combinatorics of Schubert polynomials and determinantal ideals. They proved that any subword complex is homeomorphic to a ball or a sphere and asked about their geometric realizations. 

Given a  matrix Schubert variety $ \xpi$,  it can be written as $\xpi=Y_\pi\times \mathbb{C}^q$ (where $q$ is maximal possible). Our main results are as follows. We characterize when $Y_{\pi}$ is toric (with respect to a $(\mathbb{C}^*)^{2n-1}$-action) and study the moment polytope $\Phi(\mathbb{P}(Y_\pi))$ of its projectivization. We construct regular triangulations of $\Phi(\mathbb{P}(Y_\pi))$ which we show are geometric realizations of a family of subword complexes.  Since the appearance of Knutson's and Miller's work in \cite{subword, annals} there has been a flurry of research into the geometric realization of subword complexes with progress in realizing  families of  subword complexes homeomorhpic to spheres \cite{stump,cesar,assoc,abrick,SerranoStump,subwordcluster,fans}.  The first paper which tackled  realizing a large family of subword complexes which are homeomorphic to balls is \cite{us}. Our current paper is further progress in this direction.  

The roadmap of this paper is as follows. In Section \ref{sec:schub} we define matrix Schubert varieties  $\xpi=Y_\pi\times \mathbb{C}^q$  and calculate the moment polytope $\Phi(\mathbb{P}(Y_\pi))$ of the projectivization of $Y_{\pi}$. In Section \ref{sec:polytope} we characterize when $Y_{\pi}$ is toric and construct a  regular triangulation of  $\Phi(\mathbb{P}(Y_\pi))$. In Section \ref{sec:sub} we define subword complexes and show that the aftermentioned triangulations of  $\Phi(\mathbb{P}(Y_\pi))$ are a geometric realization of a family of subword complexes homeomorphic to balls. Finally, in Section \ref{sec:deg} we show how to view the results of \cite{us} in terms of canonical triangulations of $\Phi(\mathbb{P}(Y_\pi))$.  

 
\section{Toric Schubert varieties}
\label{sec:schub}

Given a matrix Schubert variety $ \xpi$ we define a variety $Y_\pi\hookrightarrow \xpi$ such that when $Y_\pi$ is a toric variety, we can construct a regular triangulation of its corresponding polytope, which we show is a geometric realization of a family of subword complexes,  see Proposition \ref{regular} and
Theorem \ref{sc}.   Furthermore, in Theorem \ref{thm:toric} we use the diagram of  $\pi$ (defined below) to characterize the permutations $\pi$ for which $Y_\pi$ is a toric variety.

\subsection{Background and observations about matrix Schubert varieties} Let $M_n$ denote $n\times n$ matrices over $\mathbb{C}$, $B_+$ denote upper triangular invertible $n\times n$ matrices and $B_-$ denote lower triangular invertible $n\times n$ matrices. We let $\pi\in S_n$ denote both a permutation and its corresponding permutation matrix, where its $(i,j)$-th entry is 
\begin{displaymath}(\pi)_{(i,j)}=\begin{cases} 1, \text{ if } \pi(j)=i,\\ 0, \text{ else}. \end{cases}\end{displaymath}
An $n\times n$-matrix can always be reduced into a {\bf partial permutation matrix}, which is a matrix with only 0,1-entries having at most one nonzero entry on each row and column, by multiplying on the left by matrices in $B_-$ and on the right by matrices in $B_+$. The multiplication on the left corresponds to downward row operations and multiplication on the right corresponds to rightward column operations. This multiplication gives a left action of $B_-\times B_+$ on $M_n$ defined by 
\be (X,Y)\cdot M:=XMY^{-1}. \label{Bact}\ee 

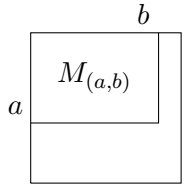
\begin{figure}[h]\centering
\begin{tikzpicture}
\draw (.85,1.4) node {$M_{(a,b)}$};
\draw (-.2,1) node {$a$};
\draw (1.5,2.25) node {$b$};
\draw (0,0)--(2,0)--(2,2)--(0,2)--(0,0);
\draw (0,.8)--(1.7,.8)--(1.7,2);
\end{tikzpicture}\caption{The submatrix $M_{(a,b)}$ of $M$.}\label{fig:submatrix}
\end{figure}

Given $1\leq a\leq m$ and $1\leq b\leq m$, let $M_{(a,b)}$ denote the upper left $a\times b$ submatrix of the matrix $M$. See Figure \ref{fig:submatrix} for an illustration. Define a {\bf rank function} of a matrix $M$ to be $r_M(a,b):=\text{rank}(M_{(a,b)})$. We then have that $M\in B_-\pi B_+$ if and only if $r_M(a,b)=r_\pi(a,b)$ for all $(a,b)\in[m]\times[m]$. 

\begin{definition} The {\bf matrix Schubert variety} of $\pi$ is $\xpi:=\overline{B_- \pi B_+}$, i.e. the Zariksi closure of its $(B_-\times B_+)$-orbit inside $M_n=\mathbb{C}^{n^2}$.
\end{definition}

Fulton studied this affine variety in \cite{MR1154177}. We summarize some of his results here. 
\begin{theorem}\cite[Proposition 3.3]{MR1154177} The matrix Schubert variety $\xpi$ is an irreducible variety of dimension $n^2-\ell(\pi)$ defined as a scheme by the equations $r_{M}(a,b)\leq r_\pi(a,b)$ for all $(a,b)\in[n]\times[n]$.
\end{theorem}

Some of these inequalities are implied by others, and Fulton described the minimal set of rank conditions.

\bd The {\bf (Rothe) diagram} of a permutation $\pi$ is the collection of boxes $D(\pi)=\{(\pi_j,i): i<j, \pi_i>\pi_j\}$. It can be visualized by considering the boxes left in the $n\times n$ grid after we cross out the boxes appearing south and east  of each 1 in the permutation matrix for $\pi$.
\ed

\begin{figure}[h]
\begin{tikzpicture}
\draw (0,0)--(2.5,0)--(2.5,2.5)--(0,2.5)--(0,0);
\draw[]
    (2.5,1.75)
 -- (0.25,1.75) node {$\bullet$}
 -- (0.25,0);
\draw[]
    (2.5,0.25)
 -- (0.75,0.25) node {$\bullet$}
 -- (0.75,0);
 \draw[]
    (2.5,0.75)
 -- (1.25,0.75) node {$\bullet$}
 -- (1.25,0); 
 \draw[]
    (2.5,2.25)
 -- (1.75,2.25) node {$\bullet$}
 -- (1.75,0);
 \draw[]
    (2.5,1.25)
 -- (2.25,1.25) node {$\bullet$}
 -- (2.25,0);
 
 \draw (0,2)--(0.5,2)--(0.5,2.5);
 \draw (0.5,2)--(1,2)--(1,2.5);
 \draw (1,2)--(1.5,2)--(1.5,2.5);
 
 \draw (1,1)--(1.5,1)--(1.5,1.5)--(1,1.5)--(1,1)--(1,0.5)--(0.5,0.5)--(0.5,1)--(1,1);
 \draw (0.5,1)--(0.5,1.5)--(1,1.5);
 
\end{tikzpicture}\caption{The diagram for $\pi=[25413]$.}
\end{figure}
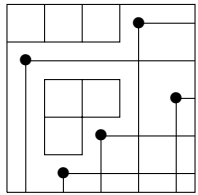

\bd {\bf Fulton's essential set}, denoted $Ess(\pi)$, is the set consisting of the south-east corners of $D(\pi)$.
\ed

\begin{theorem}\label{thm:fulton}\cite[Lemma 3.10]{MR1154177} The ideal defining the variety $\xpi$ is generated by the equations $r_M(a,b)\leq r_\pi(a,b)$ for all $(a,b)\in Ess(\pi)$.
\end{theorem}

We now define some regions inside the $(n\times n)$-grid and some varieties corresponding to these regions.
\bd The {\bf dominant piece}, denoted $dom(\pi)$, of a permutation $\pi$ is the connected component of the diagram of $\pi$ containing the box $(1,1)$, or empty if $\pi(1)=1$.
\ed 

We have that $r_\pi(a,b)=0$ if and only if $(a,b)\in dom(\pi)$. 
Therefore the dominant piece of $\pi$ consists precisely of the coordinates in $M_n$ that are 0 on $\xpi$.

\bd Let $\boldsymbol{NW(\pi)}$ denote the union over the entries north-west of some box in $D(\pi)$. Let $\boldsymbol{L(\pi)}:=NW(\pi)-dom(\pi)$ and let $\boldsymbol{L'(\pi)}:=L(\pi)-D(\pi)$.
\ed

\begin{figure}[h]
\begin{tikzpicture}
 
 \filldraw[fill=gray!20!white, draw=black] (0,0.5)--(1,0.5)--(1,1)--(1.5,1)--(1.5,2)--(0,2)--(0,0.5);
 \filldraw[fill=gray!50!white, draw=black] (0,0.5)--(0.5,0.5)--(0.5,1.5)--(1.5,1.5)--(1.5,2)--(0,2)--(0,0.5);

\draw (0,0)--(2.5,0)--(2.5,2.5)--(0,2.5)--(0,0);
\draw[]
    (2.5,1.75)
 -- (0.25,1.75) node {$\bullet$}
 -- (0.25,0);
\draw[]     (2.5,.25) -- (0.75,.25) node {$\bullet$} -- (0.75,0);
\draw[]     (2.5,1.25) -- (2.25,1.25) node {$\bullet$} -- (2.25,0);
 \draw[]     (2.5,0.75) -- (1.25,0.75) node {$\bullet$} -- (1.25,0); 
 \draw[]   (2.5,2.25) -- (1.75,2.25) node {$\bullet$} -- (1.75,0);
 
 \draw (0,2)--(0.5,2)--(0.5,2.5);
 \draw (0.5,2)--(1,2)--(1,2.5);
 \draw (1,2)--(1.5,2)--(1.5,2.5);
 
 \draw (1,1)--(1.5,1)--(1.5,1.5)--(1,1.5)--(1,1)--(1,0.5)--(0.5,0.5)--(0.5,1)--(1,1);
 
\end{tikzpicture}\caption{Given $\pi=[25413]$, $L(\pi)$ consists of all the gray boxes and $L'(\pi)$ consists of only the darker gray boxes.}\label{fig:llprime}
\end{figure}
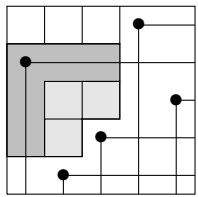

See Figure \ref{fig:llprime} for an example.

\bd Given a permutation $\pi$, let $\boldsymbol{Y_\pi}$ be the projection of $\xpi$ onto the entries inside $L(\pi)$ and let $\boldsymbol{V_\pi}$ be the projection onto the entries not north-west of any box of $D(\pi)$.
\ed 

Theorem \ref{thm:fulton} implies that the entries in $V_\pi$ are free in $\xpi$ and thus $V_\pi\cong \mathbb{C}^q$, where $q$ is the number of boxes in the region defining $V_\pi$. We have that that $\xpi=Y_\pi\times V_\pi$ and
\begin{align}\label{eq:dimyp}\begin{split}
\text{dim}(Y_\pi)&=\text{dim}(\xpi)-\text{dim}(V_\pi) \\
&= (n^2-\ell(\pi))-(n^2-|NW(\pi)|)\\
&=(n^2-|D(\pi)|)-(n^2-|NW(\pi)|) \\
&= |NW(\pi)|-D(\pi)\\
&= |L'(\pi)|.
\end{split}\end{align} 
 Consider permutations of the form $\pi=1\pi'$ where $\pi'$ is a dominant permutation of $\{2,\ldots,n\}$, i.e., the diagram of $\pi$ is a partition with north-west most box at position $(2,2)$. For these type of permutations we have that $L'(1\pi')$ is a hook.

 The $(B_-\times B^+)$-action on $\xpi$ defined in equation \eqref{Bact} restricts to a $(T^n\times T^n)$-action on $\xpi$, where $T^n$ consists of $n\times n$ diagonal matrices. This action is not faithful because scaling acts the same way on both sides, so  \begin{displaymath}\text{Stab}(T^{2n})=\{(a\cdot I,a\cdot I): a\in \mathbb{C}^*\}.\end{displaymath}
Notice that both $Y_\pi$ and $V_\pi$ are $T^{2n}$-invariant subspaces of $\xpi$. In Theorem \ref{thm:toric} we will characterize the $\pi$ for which $Y_\pi$ is a toric variety with respect to $T^{2n}/\text{Stab}(T^{2n})$ in terms of $L'(\pi)$. In other words, we will characterize the $\pi$ such that $Y_\pi$ has a dense $T^{2n}$-orbit. We will denote the quotient $T^{2n}/\text{Stab}(T^{2n})$ by $T^{2n-1}$. We note that by Fulton's realization of matrix Schubert varieties in \cite{MR1154177} as subvarieties of Scubert varieties, which are normal \cite{normal1,normal2}, it follows that matrix Schubert varieties and $Y_\pi$ are normal. 

Let $p$ be a general point of $Y_\pi$, then $\torb\subset Y_\pi$ is the affine toric variety associated to the $T^{2n}$-moment\footnote{\label{note1}The reason we use the word moment for these convex objects is because they arise in the context of symplectic and pre-symplectic geometry. For readers interested in the connection, we refer them to \cite{CdS,toricCdS,preham}.}  cone of $Y_\pi$, which we denote by $\boldsymbol{\Phi(Y_\pi)}$, and $\text{dim}(\torb)=\text{dim}(\Phi(Y_\pi))$. Since both $Y_\pi$ and $\torb$ are irreducible, if they have the same dimension they must be equal. In Theorem \ref{thm:toric} we classify when $Y_\pi$ is a toric variety by classifying the $\pi$ for which $\text{dim}(\Phi(Y_\pi))=\dim(Y_\pi)$. 

In order to describe the cone $\Phi(Y_\pi)$, we start by describing the cone $\Phi(\xpi)$ corresponding to a $T^{2n}$-orbit of a general point $q$ in $\xpi$, which we can assume without loss of generality to be $q=(1,\ldots,1)$. 
The cone $\phi(\xpi)$ is spanned by the weights by which $T^{2n}$ acts on $\torbq$. 
More explicitly, the orbit $\torbq$ is the Zariski closure of the image of a map $\varphi:T^{2n}\rightarrow \mathbb{C}^{n^2}$ where $\varphi(t)= (t^{a_{(1,1)}}q_{(1,1)},\ldots,t^{a_{(n,n)}}q_{(n,n)})$ and $\Phi(\xpi)$ is the cone spanned by the exponents $a_{(i,j)}$ of the monomials.
Notice that if $A$ and $B$ are the diagonal matrices with diagonal entries $(a_1,\ldots,a_n)$ and $(b_1,\ldots,b_n)$, respectively, then for any matrix $M$ the $(i,j)$-th entry of $AMB^{-1}$ is $a_ib_j^{-1}M_{(i,j)}$. 
We therefore have that the exponents are $x_i-y_j$, where the $x_i$ are the standard basis for $\mathbb{R}^n\times 0$, and the $y_j$ are the standard basis for $0\times \mathbb{R}^n$. In other words, the moment cone $\Phi(\xpi)$ is the cone spanned by the vectors in the set $\{x_i-y_j\mid (i,j)\in[n]\times [n]\}$.
Since $Y_\pi\hookrightarrow\xpi$ by restricting $\xpi$ to the entries inside $L(\pi)$, then $\Phi(Y_\pi)$ is the cone spanned by the set $\{x_i-y_j \mid (i,j)\in L(\pi)\}$.

Since the ideal defining $\xpi$ is homogeneous, the variety $\xpi$ is a cone, meaning that for any $z\in\xpi$ and $c\in \mathbb{C}$, we have that $cz\in\xpi$.
 We can therefore projectivize it, that is we can take 
\begin{displaymath}\mathbb{P}(\xpi):=\{[z_{(1,1)},\ldots,z_{(n,n)}] : (z_{(1,1)},\ldots,z_{(n,n)})\in \xpi\}\subset \mathbb{C}\mathbb{P}^{n^2-1}.\end{displaymath}
The same is true for $Y_\pi$. In this paper we study the {\bf moment\footnote{See footnote \ref{note1}.} polytope} $\boldsymbol{\Phi(\mathbb{P}(Y_\pi))}$ of the projectivization of $Y_\pi$. This polytope is the convex hull of $(x_i-y_j)$ for $(i,j)$ inside $L(\pi)$. The next section studies the properties of the moment polytopes $\Phi(\mathbb{P}(Y_\pi))$.

\section{Understanding the polytope $\Phi(\mathbb{P}(Y_\pi))$}
\label{sec:polytope}

This section is devoted to the study of $\Phi(\mathbb{P}(Y_\pi))={\rm ConvHull}(x_i-y_j \mid (i,j)\in L(\pi))$ for $\pi \in S_n$, the moment polytope of the projectivization of $Y_\pi$. The polytope $\Phi(\mathbb{P}(Y_\pi))$ is a root polytope, since its vertices are positive roots of type $A_{n-1}$. We set our notation for root polytopes now.

\subsection{Root polytopes and their triangulations}
A \textbf{root polytope} (of   type $A_{n-1}$) in this section will be the convex hull of  some of the points $e_i-e_j$ for $1\leq i<j \leq n$. Given a graph $G$ on the vertex set $[n]$ we associate to it the root polytope
  \be {Q}_G={\rm ConvHull}(e_i-e_j\mid (i,j) \in E(G), i<j).\ee 

We will also need a different root polytope, defined as:
\be \tilde{Q}_G={\rm ConvHull}(0, e_i-e_j\mid (i,j) \in E(G), i<j).\ee 

Note that for every $\pi \in S_n$ we have that  \textit{$L(\pi)$ is a skew Ferrers diagram}.  Given a skew Ferrers diagram $D$ with $r$ rows and $c$ columns, label its rows by $1, 2, \ldots, r$ from top to bottom and its columns by $1, 2, \ldots, c$ from left to right.  Define 
\be G_{D}=(\{x_1, \ldots, x_r, y_1, \ldots, y_c\}, \{(x_i, y_j) \mid (i,j)\in D\}).\ee  
Then \be \label{p} \Phi(\mathbb{P}(Y_\pi))=Q_{G_{L(\pi)}}.\ee
Note that an edge $(x_i, y_j) \in G_{D}$ yields the vertex $e_i-e_{r+j}$ of the root polytope $Q_{G_D}$. 

Given a drawing of a graph $G$ so that its vertices $v_1, \ldots, v_n$ are arranged in this order on a horizontal line, and its edges are drawn above this line, we say that $G$ is \textbf{noncrossing} if it has no edges $(v_i, v_k)$ and $(v_j, v_l)$ with $i<j<k<l$. A vertex $v_i$ of $G$ is said to be \textbf{nonalternating} if it has both an incoming and an outgoing edge; it is called \textbf{alternating} otherwise. The graph $G$ is alternating if all its vertices are alternating.

Since being noncrossing depends on the drawing of the graph it is essential that we  set a way to draw   $G_D$. For the purposes of this paper the vertices of $G_D$  are drawn from  left to right in the following order: $x_r, \ldots, x_1, y_c, \ldots, y_1$.

\bl \label{ggd-tri} Given a skew diagram $D$ for which $G_D$ has $k$ components, the root polytope ${Q}_{G_{D}}=\bigcup_{F}{Q}_{F}$, where the union runs over all noncrossing alternating spanning forests of $G_{D}$ with $|V(G_D)|-k$ edges and the simplices ${Q}_{F}$ are interior disjoint and of the same dimension as ${Q}_{G_{D}}$.   
\el

Before proving Lemma \ref{ggd-tri},  we prove  Lemma  \ref{bij}, which gives another way to view the triangulation described in Lemma \ref{ggd-tri}. Special cases of these lemmas were  used in \cite{greene}, and  the triangulation of the polytope $Q_{G_D}$ into simplices labeled by the lattice paths appearing in Lemma \ref{bij} was proven by different methods in \cite{reiner}.

\bl \label{bij} Given a skew diagram $D$ with  $k$ components $D_1, \ldots, D_k$, the noncrossing alternating spanning forests of $G_D$  with $|V(G_D)|-k$ edges are in bijection with $k$-tuples of lattice paths   $(p_1, \ldots, p_k)$, where $p_i$ is a lattice path inside $D_i$ starting at  $D_i$'s southwestmost box and  ending at at  $D_i$'s northeastmost box  taking steps either one unit south or west.
\el

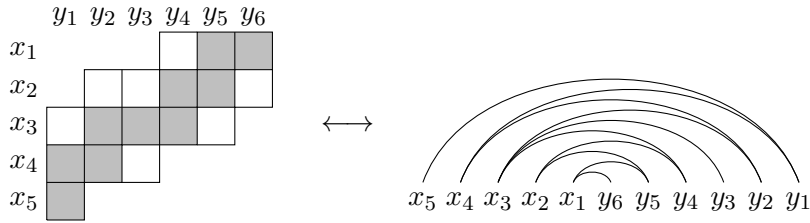
\begin{figure}[h]
 \centering
 \begin{tikzpicture}
  \filldraw[fill=gray!50!white, draw=black] (0,0)--(0,1)--(.5,1)--(.5,1.5)--(1.5,1.5)--(1.5,2)--(2,2)--(2,2.5)--(3,2.5)--(3,2)--(2.5,2)--(2.5,1.5)--(2,1.5)--(2,1)--(1,1)--(1,0.5)--(0.5,0.5)--(0.5,0)--(0,0);
  
  \draw[]     (0,.5) -- (.5,.5)--(.5,1)--(1,1)--(1,2);
    \draw[]     (1,.5) -- (1.5,.5) -- (1.5,1.5);
      \draw[]     (0,1) -- (0,1.5)--(.5,1.5)--(.5,2)--(1.5,2)-- (1.5,2.5)--(2,2.5);
      \draw[]     (1.5,1.5) -- (2,1.5)--(2,2)--(2.5,2)--(2.5,2.5);
      \draw[]     (3,2) -- (3,1.5)--(2.5,1.5)--(2.5,1)--(2,1);
      
      \draw (-.3,2.25) node {$x_1$};
      \draw (-.3,1.75) node {$x_2$};
      \draw (-.3,1.25) node {$x_3$};
      \draw (-.3,.75) node {$x_4$};
      \draw (-.3,.25) node {$x_5$};
      
      \draw (2.75,2.7) node {$y_6$};      
      \draw (2.25,2.7) node {$y_5$};
      \draw (1.75,2.7) node {$y_4$};
      \draw (1.25,2.7) node {$y_3$};
      \draw (.75,2.7) node {$y_2$};
      \draw (.25,2.7) node {$y_1$};
      
      \draw (4,1.25) node {$\longleftrightarrow$};
      
      \path (5,0.5) edge [bend left=70] (10,.5);
\path (5.5,0.5) edge [bend left=70] (10,0.5);
\path (5.5,0.5) edge [bend left=70] (9.5,0.5);
\path (6,0.5) edge [bend left=70] (9.5,0.5);
\path (6,0.5) edge [bend left=70] (9,0.5);
\path (6,0.5) edge [bend left=70] (8.5,0.5);
\path (6.5,0.5) edge [bend left=70] (8.5,0.5);
\path (6.5,0.5) edge [bend left=70] (8,0.5);
\path (7,0.5) edge [bend left=70] (7.5,0.5);
\path (7,0.5) edge [bend left=70] (8,0.5);

\draw[] (5,.25) node {$x_5$};
\draw[] (5.5,.25) node {$x_4$};
\draw[] (6,.25) node {$x_3$};
\draw[] (6.5,.25) node {$x_2$};
\draw[] (7,.25) node {$x_1$};

\draw[] (7.5,.25) node {$y_6$};
\draw[] (8,.25) node {$y_5$};
\draw[] (8.5,.25) node {$y_4$};
\draw[] (9,.25) node {$y_3$};
\draw[] (9.5,.25) node {$y_2$};
\draw[] (10,.25) node {$y_1$};
                  
 \end{tikzpicture} \caption{The correspondence between noncrossing alternating spanning trees of $G_D$ and  lattice paths    from $(1,c)$ to $(r,1)$ inside $D$ that take steps either one unit south or west. }
 \label{fig:bij}
 \end{figure}

 \proof The bijection is given by the map that takes a noncrossing alternating spanning forest $F=(\{x_r, \ldots, x_1, y_c, \ldots, y_1\}, \{(x_i, y_j) \mid (i,j)\in S(F)\})$ of $G_D$ to the $k$-tuple of paths $(p_1, \ldots, p_k)$, where the set of edges of $p_1\cup \cdots\cup p_k$ is given by the set $S(F)$. See Figure \ref{fig:bij} for an example. \qed

\noindent \textit{Proof of Lemma \ref{ggd-tri}.} If $D$ is connected, then 
$G_D$ is a connected bipartite graph on the parts $\{x_1, \ldots, x_r\}$ and $\{y_1, \ldots, y_c\}$. Thus, $G_D \subset K_{r,c}$, the bipartite graph with parts of sizes $r$ and $c$. By \cite[Proposition 13.3]{p1} the polytopes $Q_{K_{r,c}}$ are the facets of $\tilde{Q}_{K_{r+c}}$, where $K_{r+c}$ is the complete graph on $r+c$ vertices, and since $G_D$ is connected we have that $Q_{G_D}$ and  $Q_{K_{r,c}}$ are of the same dimension.  Thus, it suffices to prove that $\tilde{Q}_{G_{D}}=\bigcup_{T}\tilde{Q}_{T}$, where the union runs over all noncrossing alternating trees of $G_{D}$ and the simplices $\tilde{Q}_{T}$ are interior disjoint and of the same dimension as $\tilde{Q}_{G_{D}}$.  The latter statement follows  
by application of \cite[Lemma 12.6]{p1}.  

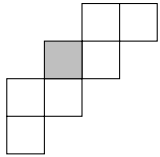
\begin{figure}[h]
\begin{tikzpicture}
  \draw (0,0.5)--(1,0.5)--(1,1)--(0,1)--(0,0)--(0.5,0)--(0.5,1);
  \draw(1,1.5)--(2,1.5)--(2,2)--(1,2)--(1,1)--(1.5,1)--(1.5,2);
  \filldraw[fill=gray!50!white, draw=black] (0.5,1)--(0.5,1.5)--(1,1.5)--(1,1)--(0.5,1); 
\end{tikzpicture}
\caption{For $D$ the unshaded region, $G_D$ is disconnected. In this case, $D'$ equals $D$ together with the shaded square.}\label{fig:D'}
\end{figure}

If $D$ is not connected, there exist partitions $\lambda, \mu, \omega$ such that $D=\lambda\setminus \mu$ and $D':=\lambda\setminus\omega$ is the unique smallest such  that $G_{D'}$ is connected; see Figure \ref{fig:D'}.  The number of rows and columns of $D$ and $D'$ are the same, $r$ and $c$, respectively,  and $G_D \subset G_{D'}$.  Since  $G_{D'}$ is connected the first paragraph of this proof applies to it.   Denoting by $B$ the set of boxes in $D'$ that are not in $D$, it follows that every noncrossing alternating spanning tree of $G_{D'}$ contains the set of  edges $E_B:=\{(x_i,y_j) \mid (i,j) \in B\}$. An easy way to see this is to notice that every lattice path described in Lemma \ref{bij} must contain the boxes in $B$. Note that if we delete the edges in $E_B$ from the noncrossing alternating spanning trees of $G_{D'}$ we exactly obtain the forests described in the lemma. 
We get that in the triangulation of $Q_{G_{D'}}$ proven in the first paragraph  the vertices corresponding to the edges in $E_B$ are all cone points, and a triangulation of $Q_{G_{D'}}$  can be obtained from a triangulation of $Q_{G_{D}}$ by adding the vertices corresponding to the edges in $E_B$ (and thereby raising the dimension by $|E_B|$).  It is clear that  $\bigcup_{F}{Q}_{F} \subset  {Q}_{G_{D}}$ and the simplices ${Q}_{F}$ are interior disjoint, intersect on common faces  and  are of the same dimension as ${Q}_{G_{D}}$, since they are the restrictions of the noncrossing alternating trees of $G_{D'}$ to $G_D$.   It also follows that ${Q}_{G_{D}} \subset \bigcup_{F}{Q}_{F},$ since if a triangulation of ${Q}_{G_{D}}$ were to include simplices other than the ${Q}_{F}$'s, then the triangulation of ${Q}_{G_{D}}$ would include  simplices other than the ${Q}_{T}$'s, which would contradict the statement proven in the first paragraph. \qed

\medskip

We call the triangulation of ${Q}_{G_{D}}$ given in Lemma \ref{ggd-tri} the \textbf{noncrossing alternating triangulation}, or {\bf NAT} for short.  Recall that a triangulation of a polytope $P$ is {\bf regular} if there exists a concave  piecewise linear function $f:P \rightarrow \mathbb{R}$ such that the regions of linearity of $f$ are the maximal simplices in the triangulation.   

\bp \label{regular} For a skew diagram $D$, the NAT triangulation of $Q_{G_D}$ described in Lemma \ref{ggd-tri} is a regular triangulation.\ep

\proof We start by describing how to extend a function defined on a set $A$ of points into a piecewise linear function defined on their convex hull.
Let $f:A \rightarrow \mathbb{R}$ be a function on the set $A$ and consider the polytope $P={\rm ConvHull}(A)$. 
Let $\overline{P}={\rm ConvHull}((a, f(a)) \mid a \in A)$ and define then $f(p)={\rm max}\{x \mid (p, x) \in \overline{P}\}$, $p \in P$.   The  function $f:P \rightarrow \mathbb{R}$ is concave by definition. 

Consider the root polytope $Q_{G_D}$ with  vertices $e_i-e_{r+j}$, for  $(i, j) \in D$. 
Let  $f(e_i-e_{r+j})=(i-(r+j))^2$ for $(i, j) \in D$.  Extend this to a concave piecewise linear function on $Q_{G_D}$ as described above. It can be shown that this functions yields the NAT triangulation of $Q_{G_D}$, thus showing its regularity.
\qed

\subsection{Characterizing when $Y_{\pi}$ is a toric variety} Now we are ready to use the above lemmas in order to characterize when $Y_{\pi}$ is a toric variety.

\bl \label{cor-dim} Given a skew diagram $D$ with $r$ rows and $c$ columns, for which $G_D$ has $k$ components, the dimension of ${Q}_{G_{D}}$ is $r+c-k-1$. 
\el 

\proof  The dimension of ${Q}_{F}$ for some noncrossing alternating spanning forest $F$ of $G_{D}$ is $r+c-k-1$, since $F$ has $r+c-k$ edges,  which correspond to the vertices of a $(r+c-k-1)$-dimensional simplex.
  Since  by Lemma \ref{ggd-tri}, the dimension of ${Q}_{G_{D}}$ equals the dimension of ${Q}_{F}$ for some noncrossing alternating spanning forest $F$ of $G_{D}$ with $|V(G_D)|-k$ edges, the statement follows. \qed

\bt \label{thm:toric} $Y_\pi$ is a toric variety with respect to the $T^{2n-1}$-action if and only if $L'(\pi)$ consists of disjoint hooks that do not share a row or a column with each other.
\et

\proof In Equation (\ref{eq:dimyp}) we computed the dimension of $Y_\pi$. Lemma \ref{cor-dim}    yields that  the dimension of $\Phi(\mathbb{P}(Y_\pi))$  equals $|L'(\pi)|-1$ if and only if $L'(\pi)$ consists of disjoint hooks that do not share a row or a column with each other. This suffices to prove the theorem.
\qed
\medskip

An immediate corollary of Theorem \ref{thm:toric} is the following.

\bc If $\pi'$ is a dominant permutation on $2,3,\ldots,n$ then $Y_{1\pi'}$ is a toric variety.
\ec


\section{On geometric realizations of subword complexes}
\label{sec:sub}
In this section we show that the NAT triangulations studied in the previous section  geometrically realize certain subword complexes. Thus,  we extend the family of subword complexes homeomorphic to balls that have known geometric realizations as triangulations of polytopes.

The symmetric group $S_n$ is generated by the adjacent transpositions $s_1,\ldots,s_{n-1}$, where $s_i$ transposes $i \leftrightarrow i+1$. Let $Q=(q_1,\ldots,q_\sizeQ)$ be a word in $\{s_1,\ldots,s_{n-1}\}$, i.e., $Q$ is an ordered sequence. A {\bf subword} $J=(r_1,\ldots,r_\sizeQ)$ of $Q$ is a word obtained from $Q$ by replacing some of its letters by $-$. There are a total of $2^{|Q|}$ subwords of $Q$. Given a subword $J$, we denote by $Q\setminus J$ the subword with $k$-th entry equal to $-$ if $r_k\neq -$ and equal to $q_k$ otherwise for, $k=1,\ldots,\sizeQ$. For example, $J=(s_1,-,s_3,-,s_2)$ is a subword of $Q=(s_1,s_2,s_3,s_1,s_2)$ and $Q\setminus J=(-,s_2,-,s_1,-)$. Given a subword $J$ we denote by $\prod J$ the product of the letters in $J$, from left to right, with $-$ behaving as the identity.

\bd \cite{subword,annals}
 Let $Q=(q_1,\ldots,q_\sizeQ)$ be a word in $\{s_1,\ldots,s_{n-1}\}$ and $\pi\in S_n$. The \textbf{subword complex} $\Delta(Q,\pi)$ is the simplicial complex on the vertex set $Q$ whose facets are the subwords $F$ of $Q$ such that the product $\prod(Q\setminus F)$ is a reduced expression for $\pi$.
\ed

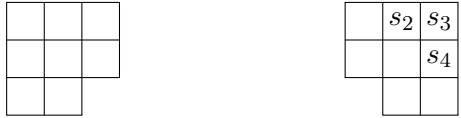
\begin{figure}[h]
\begin{tikzpicture}

 \draw (1,1.5)--(1,0)--(0,0)--(0,1.5)--(1.5,1.5)--(1.5,.5)--(0,.5);
 
 \draw[] (.5,1.5)--(.5,0);
 \draw[] (0,1)--(1.5,1);

\begin{scope}[xscale=-1,xshift=-9cm]
 \draw (4,1.5)--(4,0)--(3,0)--(3,1.5)--(4.5,1.5)--(4.5,.5)--(3,.5);
 
 \draw[] (3.5,1.5)--(3.5,0);
 \draw[] (3,1)--(4.5,1);
 \end{scope}
 
\node at (5.25,1.25) {$s_2$};
\node at (5.75,1.25) {$s_3$};
\node at (5.75,.75) {$s_4$};
 
\end{tikzpicture}\caption{On the left we have $L(\pi)$ and on the right $\overline{L(\pi)}$ for $\pi=[14523]$. Note that $p(\pi)=s_4s_2s_3=[13524]$.}\label{fig:mirrorL}
\end{figure}

Given a permutation $\pi \in S_n$, let $\overline{L(\pi)}$ be the mirror image of  the skew shape $L(\pi)$. 
 Let $Q(\L)$ be the word given by reading the transpositions in the boxes of $\L$ from left to right and from bottom to top. Let $P(\pi)=\L-B(\pi)$ where $B(\pi)$ is as follows. In each connected part of $\L$ draw the lowestmost path from its top left box to its bottommost rightmost box. These boxes constitute $B(\pi)$. Let $p(\pi)$ be the permutation obtained from reading the transpositions in the boxes of $P(\pi)$ from left to right and from bottom to top. See Figure \ref{fig:mirrorL}.

\bt \label{sc}For $\pi \in S_n$, the noncrossing alternating triangulation of  $Q_{G_{L(\pi)}}$ is a geometric realization of the subword complex $\Delta(Q(\L), p(\pi))$. 
\et



\proof The dimension of $\Delta(Q(\L), p(\pi))$ is $|B(\pi)|-1$ by definition. On the other hand the dimension of  $Q_{G_{L(\pi)}}$ equals this same thing by Lemma \ref{ggd-tri} (this is clear through  Lemma \ref{bij}). Viewing the noncrossing alternating spanning forests labeling the top dimensional simplices in NAT as a union  of lattice paths via  Lemma \ref{bij}, we get  a bijection  to reduced pipe dreams of $p(\pi))$ by filling by elbow the boxes that are in the paths labeling the simplices of NAT  (or rather in their mirror images in $\L$) and putting crosses in the rest of $\L$. Since $\Delta(Q(\L), p(\pi))$ and the NAT of $Q_{G_{L(\pi)}}$  are pure simplicial complexes, the above proves the theorem.
\qed

\section{Degeneration of moment polytopes into acyclic root polytopes}
\label{sec:deg}

In \cite[Theorem 1]{us} we have shown that the pipe dream complex $PD(\pi)$ of a  permutations $\pi=1\pi'$, with  $\pi'$ dominant, can be geometrically realized as canonical triangulations of acyclic root polytopes $\mathcal{P}(T(\pi))$.  The \textbf{pipe dream complex} $PD(\pi)$ is the subword complex $\Delta(Q,\pi)$ corresponding to the triangular word $Q=(s_{n-1},s_{n-2},s_{n-1},\ldots,s_1,s_2,\ldots,s_{n-1})$ and $\pi$.  In this section we review   \cite[Theorem 1]{us} and relate the moment polytopes $\Phi(\mathbb{P}(Y_\pi))$ to the acyclic root polytopes $\mathcal{P}(T(\pi))$. In turn, we give another way to see $PD(\pi)$ as a triangulation of $\Phi(\mathbb{P}(Y_\pi))$.

\subsection{Acyclic root polytopes} Let $G$ be an acyclic   graph on the vertex set $[n+1]$.  Define $$\mathcal{V}_G=\{e_i-e_j \mid  (i, j) \in E(G), i<j\}, \mbox{ a set of vectors associated to $G$;}$$

 $${\rm cone}(G)=\langle \mathcal{V}_G \rangle :=\{\sum_{ e_i-e_j \in \mathcal{V}_G}c_{ij} (e_i-e_j) \mid  c_{ij}\geq 0\}, \mbox{ the cone associated to $G$; and } $$  
  $$\overline{\mathcal{V}}_G=\Phi^+ \cap {\rm cone}(G), \mbox{ all the positive roots of type $A_n$ contained in ${\rm cone}(G)$}, $$
   where $\Phi^+=\{e_i-e_j  \mid1\leq i<j \leq n+1\}$ is the set of         positive roots of type $A_n$.  
    
The \textbf{acyclic root polytope} $\mathcal{P}(G)$ associated to the acyclic graph $G$ is  
  
 \begin{equation} \label{eq11} \mathcal{P}(G)=\textrm{ConvHull}(0, e_i-e_j \mid e_i-e_j \in \overline{\mathcal{V}}_G).\end{equation}

\bt \label{?} \cite{root1} Let $T_1, \ldots, T_k$ be the noncrossing alternating spanning trees of the directed transitive closure of the acyclic graph $G$. Then $\P(T_1), \ldots, \P(T_k)$ are top dimensional simplices in a regular  triangulation of $\P(G)$ called the \textbf{canonical triangulation}.
\et
 
 The main tool developed in \cite{root1} which is used to construct the canonical triangulation of Theorem \ref{?} is the subdivision algebra. Subdivision algebras have since been utilized in solving various problems in \cite{h-poly2, h-poly1, greene, pipe, prod, mm}.

 \bt  \label{thm:intro} \cite[Theorem 1]{us} Let $\pi=1\pi' \in S_n$, where $\pi'$ is dominant. Let $\mathcal{C}^2(\pi)$ be the core of $PD(\pi)$ coned over twice. The   canonical triangulation of the root polytope $\P(T(\pi))$   is a geometric realization of $\mathcal{C}^2(\pi)$.
 \et

\subsection{Definition of $T(\pi)$ (appearing in Theorem \ref{thm:intro})} We now review the original description from \cite{us} of $T(\pi)$ used in Theorem \ref{thm:intro}.  For a definition of a reduced pipe dream, also called RC-graphs, as well as further background on it see \cite{rc, fom-kir, subword} or \cite{us}. 

\begin{figure}[h]
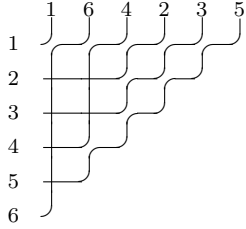
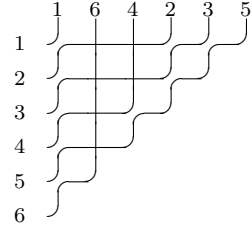

    \centering
    \begin{subfigure}[b]{0.4\textwidth}
        \centering
      $
\begin{array}{ccccccc}
       &\perm1{}&\perm6{}&\perm4{}&\perm2{}&\perm3{}&\perm5{}\\
\petit1&   \jr   &   \jr  &   \jr  &   \jr  &   \jr   &  \je   \\
\petit2&   \+  &   \+ &   \jr &   \jr    &   \je  &\\
\petit3&   \+  &   \+ &   \jr   &    \je    &\\
\petit4&   \+  &    \jr    &  \je      &\\
\petit5&   \+  &     \je   &        &\\
\petit6&   \je  &        &        &\\
\end{array}
$
        \caption{Aligned left is the bottom reduced pipe dream}
    \end{subfigure}
    \hfill
    \begin{subfigure}[b]{0.4\textwidth}
        \centering      $
\begin{array}{ccccccc}
       &\perm1{}&\perm6{}&\perm4{}&\perm2{}&\perm3{}&\perm5{}\\
\petit1&   \jr   &   \+  &   \+  &   \jr  &   \jr   &  \je   \\
\petit2&   \jr  &   \+ &   \+ &   \jr    &   \je  &\\
\petit3&   \jr  &   \+ &   \jr   &    \je    &\\
\petit4&   \jr  &    \+    &  \je      &\\
\petit5&   \jr  &     \je   &        &\\
\petit6&   \je  &        &        &\\
\end{array}
$
        \caption{Aligned up is the top reduced pipe dream}

    \end{subfigure}
    \hfill
    \caption{Two reduced pipe dreams for [164235] obtained by aligning the diagram to the left and to the top.}        \label{fig:topbottomrpd}
\end{figure}

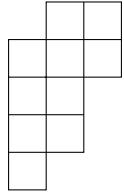
\begin{figure}[h]
\begin{tikzpicture}
 
 \draw (0,1.5)--(1.5,1.5)--(1.5,2.5)--(0.5,2.5)--(0.5,0)--(0,0)--(0,2)--(1.5,2);
 \draw (0,0.5)--(1,0.5)--(1,2.5);
 \draw (0,1)--(1,1);

\end{tikzpicture}\caption{The core region of $[164235]$}\label{fig:crpi}
\end{figure}

Since we are only considering permutations of the form $1\pi'$ with $\pi'$ a dominant permutation, core$(1\pi')$ is easy to describe. Given a diagram of a permutation there are two natural reduced pipe dreams for $\pi$, referred to as the \textbf{bottom reduced pipe dream of $\pi$} and the  \textbf{top reduced pipe dream of $\pi$}, one obtained by aligning the diagram to the left and replacing the boxes with crosses and the other one by aligning the diagram up. See Figure \ref{fig:topbottomrpd}. The core of $1\pi'$ is the simplicial complex obtained by restricting $PD(\pi)$ to the vertices corresponding to the positions of the crosses in the superimposition of these two pipe dreams. We refer to the region itself as the \textbf{core region}, and denote it by cr$(\pi)$. See Figure \ref{fig:crpi} for an example.

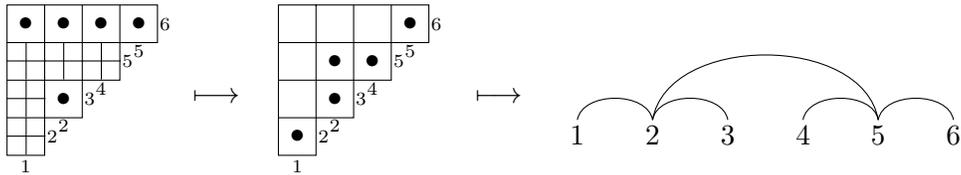
\begin{figure}[h]
    \centering
    \begin{subfigure}[h]{2.5cm}
        
          \begin{tikzpicture}
 
 \draw (0.5,2)--(0.5,0)--(0,0)--(0,2)--(2,2)--(2,1.5)--(0,1.5);
 \draw (0,0.5)--(1,0.5)--(1,2);
 \draw (0,1)--(1.5,1)--(1.5,2);

\draw (.25,1.75) node {$\bullet$};
\draw (.75,.75) node {$\bullet$};
\draw (.75,1.75) node {$\bullet$};
\draw (1.25,1.75) node {$\bullet$};
\draw (1.75,1.75) node {$\bullet$};
 
 \draw (0.26,0.15) node {$\textcross$};
 \draw (0.26,0.65) node {$\textcross$};
  \draw (0.26,1.15) node {$\textcross$};
   \draw (0.76,1.15) node {$\textcross$};
    \draw (1.26,1.15) node {$\textcross$};

\draw (0.25,-0.15) node {{\tiny 1}};
 \draw (.6,0.25) node {{\tiny 2}};
 \draw (.75,0.39) node {{\tiny 2}};
 \draw (1.1,0.75) node {{\tiny 3}};
 \draw (1.25,0.89) node {{\tiny 4}};
  \draw (1.6,1.25) node {{\tiny 5}};
 \draw (1.75,1.39) node {{\tiny 5}};
 \draw (2.1,1.75) node {{\tiny 6}};
 
\end{tikzpicture}

    \end{subfigure}
    \begin{subfigure}[h]{1cm}
        $\displaystyle{\lmp{}}$
        
            \end{subfigure}
        \begin{subfigure}[h]{2.5cm}
        
          \begin{tikzpicture}
 
 \draw (0.5,2)--(0.5,0)--(0,0)--(0,2)--(2,2)--(2,1.5)--(0,1.5);
 \draw (0,0.5)--(1,0.5)--(1,2);
 \draw (0,1)--(1.5,1)--(1.5,2);

\draw (.25,.25) node {$\bullet$};
\draw (.75,.75) node {$\bullet$};
\draw (.75,1.25) node {$\bullet$};
\draw (1.25,1.25) node {$\bullet$};
\draw (1.75,1.75) node {$\bullet$};

\draw (0.25,-0.15) node {{\tiny 1}};
 \draw (.6,0.25) node {{\tiny 2}};
 \draw (.75,0.39) node {{\tiny 2}};
 \draw (1.1,0.75) node {{\tiny 3}};
 \draw (1.25,0.89) node {{\tiny 4}};
  \draw (1.6,1.25) node {{\tiny 5}};
 \draw (1.75,1.39) node {{\tiny 5}};
 \draw (2.1,1.75) node {{\tiny 6}};
 
\end{tikzpicture}

    \end{subfigure}
    \begin{subfigure}[h]{1cm}
        $\displaystyle{\lmp{}}$
        
            \end{subfigure}
    \begin{subfigure}[h]{5.4cm}

        \begin{tikzpicture}
 
\path (1,0) edge [bend left=80] (2,0);
\path (2,0) edge [bend left=80] (3,0);
\path (2,0) edge [bend left=80] (5,0);
\path (4,0) edge [bend left=80] (5,0);
\path (5,0) edge [bend left=80] (6,0);

\draw[] (1,-0.2) node {1};
\draw[] (2,-0.2) node {2};
\draw[] (3,-0.2) node {3};
\draw[] (4,-0.2) node {4};
\draw[] (5,-0.2) node {5};
\draw[] (6,-0.2) node {6};
 
\end{tikzpicture}
        
    \end{subfigure}
\caption{Let $\pi=[15342]$, on the left we see the bottom reduced pipe dream for $\pi$ drawn inside $R(\pi)$ with dots instead of elbows, this gives the labeling of the boundary. We then drop the dots to the south to get the dots encoding $T(\pi)$, which is depicted on the right.}\label{fig:cvs}\end{figure}

Let $\pi=1\pi'$, where $\pi'$ is dominant. Denote by $\S(\pi)$ the subword complex which is the core$(\pi)$ coned over the vertex of $PD(\pi)$ corresponding to the entry $(1,1)$. Denote the region which is the union of $(1,1)$ and cr$(\pi)$ by $R(\pi)$. In order to determine the tree $T(\pi)$, we will label the southeast boundary with numbers and we will place some dots in $R(\pi)$, see Figure \ref{fig:cvs}. The boundary of the core region starting from the southwest (SW) corner of it to the northeast (NE) corner can be described as a series of east (E) and north (N) steps. Let $A$ be the set consisting of all the $N$ steps together with some $E$ steps. The step $E\in A$ if the bottom reduced pipe dream is bounded by $E$ but not by the $N$ step directly preceding $E$. As we traverse this lower boundary from the SW corner we write the numbers $1,\ldots, m$ in increasing fashion below the $E$ steps and to the right of the $N$ steps that belong to $A$. For the $E$ steps that we did not assign a number, we consider their number to be the number assigned to the $N$ step directly preceding them. Consider the bottom reduced pipe dream drawn inside $R(\pi)$ and with elbows replaced by dots. Drop these dots south. Define $T(\pi)$ to be the tree on $m$ vertices such that there is an edge between vertices $i<j$ if there is a dot in the entry in the column of the E step labeled $i$ and in the row of the N step labeled $j$.

 \subsection{Alternative description of  $T(\pi)$ and degeneration of $\mpy$ to $\P(T(\pi))$} Here we give an alternative way to describe the above through the following lemmas. 

\bp\label{prop:cress} Let $\pi=1\pi'$ with $\pi'$ dominant. Then the union of $(1,1)$ and cr$(\pi)$ consists of all entries north-west of some box in $D(\pi)$ that are not inside Fulton's essential set. In equations, this says that $R(\pi)=NW(\pi)-Ess(\pi)$.
\ep

\proof It is clear that $R(\pi)\subset NW(\pi)$. Given $(i,j)\in Ess(\pi)$, we have that neither the bottom or the top reduced pipe dreams of $\pi$ have a cross on this position since $(i,j)$ is a south-east corner of $D(\pi)$ and so $(i,j)\notin R(\pi)$. Since $D(\pi)$ is a partition, then $D(\pi)$ restricted to a row or to a column is connected. Therefore, any entry strictly west of an entry in $Ess(\pi)$ must be in the bottom reduced pipe dream of $\pi$ and thus inside $R(\pi)$. Similarly, any entry strictly north of $Ess(\pi)$ must be in the top reduced pipe dream.
\qed

Recall that the graph associated to cr$(\pi)$ for $\pi=1\pi'$ with $\pi'$ dominant has as vertices the elements of the set $A$ which consists of some $E$ and $N$ steps that bound the $SW$ boundary of cr$(\pi)$.

\bp\label{prop:esslab} Let $\pi=1\pi'$ with $\pi'$ dominant. An $E$ step is in $A$ if and only if it does not bound a box in $Ess(\pi)$.
\ep

\proof Proposition \ref{prop:cress} tells us that given an $E$ step the $N$ step directly preceding $E$ bounds $R(\pi)$ if and only if $E$ bounds a box in $Ess(\pi)$.
\qed

\medskip

\begin{figure}[h]
\begin{tikzpicture}
\draw (0,0)--(2.5,0)--(2.5,2.5)--(0,2.5)--(0,0);
\draw[]
    (2.5,2.25)
 -- (0.25,2.25) node {$\bullet$}
 -- (0.25,0);
 \draw[]
    (2.5,.75)
 -- (0.75,.75) node {$\bullet$}
 -- (0.75,0);
 \draw[]
    (2.5,.25)
 -- (1.25,.25) node {$\bullet$}
 -- (1.25,0);
 \draw[]
    (2.5,2)
 -- (1.75,2) node {$\bullet$}
 -- (1.75,0);
  \draw[]
    (2.5,1.5)
 -- (2.25,1.5) node {$\bullet$}
 -- (2.25,0);
 
 \filldraw[fill=gray!20!white, draw=black] (.5,1)--(1.5,1)--(1.5,2)--(.5,2)--(.5,1);
 
 \draw[] (1,1)--(1,2);
 \draw[] (.5,1.5)--(1.5,1.5);
 
\end{tikzpicture}\caption{The diagram for $\pi=[14523]$.}\label{fig:newtdig}
\end{figure}
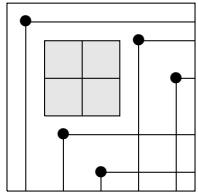

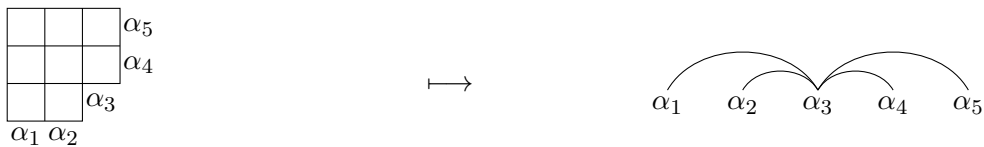
\begin{figure}[h]   
 \centering
    \begin{subfigure}[h]{0.4\textwidth}
        \centering
\begin{tikzpicture}
 
 \draw (0,1)--(1.5,1)--(1.5,2)--(0,2)--(0,.5)--(1,.5)--(1,2);
 \draw (0.5,2)--(.5,0.5);
 \draw (0,1.5)--(1.5,1.5);
 
 \draw (0.25,0.3) node {$\alpha_1$};
  \draw (0.75,0.3) node {$\alpha_2$};
   \draw (1.25,0.75) node {$\alpha_3$};
    \draw (1.75,1.25) node {$\alpha_4$};
     \draw (1.75,1.75) node {$\alpha_5$};
 
\end{tikzpicture}
 \end{subfigure}
    \hfill
    \begin{subfigure}[h]{0.1\textwidth}
        \centering   
        $\displaystyle{\lmp{}}$
        
            \end{subfigure}
    \hfill
    \begin{subfigure}[h]{0.4\textwidth}
        \centering
 
        \begin{tikzpicture}
 
\path (1,0) edge [bend left=60] (3,0);
\path (2,0) edge [bend left=60] (3,0);
\path (3,0) edge [bend left=60] (4,0);
\path (3,0) edge [bend left=60] (5,0);

\draw[] (1,-0.2) node {$\alpha_1$};
\draw[] (2,-0.2) node {$\alpha_2$};
\draw[] (3,-0.2) node {$\alpha_3$};
\draw[] (4,-0.2) node {$\alpha_4$};
\draw[] (5,-0.2) node {$\alpha_5$};
 
\end{tikzpicture}

    \end{subfigure}
    \hfill          
\caption{On the left we have $NW(\pi)-Ess(\pi)$ for $\pi=[14523]$ with its SW boundary labelled by $A$. On the right we have the corresponding tree $T(\pi)$.}\label{fig:newtmap}
\end{figure}

These propositions allow us to rephrase the construction of $T(\pi)$, for permutations $\pi=1\pi'$, $\pi'$ dominant. 
We construct the graph $T(\pi)$ by looking at the $E$ and $N$ steps bounding the $SE$ boundary of $NW(\pi)-Ess(\pi)$. Let $A$ be the set consisting of all the $N$ steps together with the $E$ steps that do not bound a box in $Ess(\pi)$. Suppose $|A|=m$, as we transverse this lower boundary from the $SW$ corner, we write $\alpha_1,\ldots,\alpha_m$ in increasing fashion below the $E$ steps and to the right of the $N$ steps that belong to $A$. For the $E$ steps that we did not assign an $\alpha_i$, we consider their label to be the $\alpha_i$ assigned to the $N$ step directly preceding them. Let $\overline{T}(\pi)$ be the tree with vertices $V=\{\alpha_1,\ldots,\alpha_m\}$ and edges corresponding to the entries in $NW(\pi)-Ess(\pi)$. More precisely $(\alpha_i,\alpha_j)$ is an edge of $\overline{T}(\pi)$ if the entry $(i,j)\in NW(\pi)-Ess(\pi)$ is north of the $E$ step labeled $\alpha_j$ and east of the $N$ step labeled $\alpha_i$. See Figures \ref{fig:newtdig} and \ref{fig:newtmap} for an example.

Note that Propositions \ref{prop:cress} and \ref{prop:esslab} tell us that for $\pi=1\pi'$ with $\pi'$ dominant this construction coincides with the one we had before.

\bt\label{thm:tortoroot} Given $\pi=1\pi'$,  with $\pi'$ dominant, the moment polytope $\mpy$ can be degenerated into the root polytope $\P(T(\pi))$.
\et

\begin{figure}[h]
    \centering
    \begin{subfigure}[h]{0.4\textwidth}
        \centering
\begin{tikzpicture}
 
 \draw (0,1)--(4.5,1)--(4.5,2)--(0,2)--(0,.5)--(3,.5)--(3,2);
 \draw (1.5,2)--(1.5,0.5);
 \draw (0,1.5)--(4.5,1.5);
 \draw (3,0.5)--(4.5,0.5)--(4.5,1);
   \filldraw[fill=gray!20!white, draw=black]  (3,0.5)--(4.5,0.5)--(4.5,1)--(3,1)--(3,0.5);

 \draw (0.75,0.75) node {$x_3-y_1$};
 \draw (0.75,1.25) node {$x_2-y_1$};
 \draw (0.75,1.75) node {$x_1-y_1$};
 \draw (2.25,0.75) node {$x_3-y_2$};
 \draw (2.25,1.25) node {$x_2-y_2$};
 \draw (2.25,1.75) node {$x_1-y_2$};
 \draw (3.75,0.75) node {$x_3-y_3$}; 
 \draw (3.75,1.25) node {$x_2-y_3$};
 \draw (3.75,1.75) node {$x_1-y_3$};

\end{tikzpicture}
 \end{subfigure}
    \hfill
    \begin{subfigure}[h]{0.1\textwidth}
        \centering   
        $\displaystyle{\lmp{K}}$
        
            \end{subfigure}
    \hfill
    \begin{subfigure}[h]{0.4\textwidth}
        \centering    
        
        \begin{tikzpicture}
 
      \filldraw[fill=gray!20!white, draw=black] (3,0.5)--(4.5,0.5)--(4.5,2)--(3,2)--(3,0.5);

 \draw (0,1)--(4.5,1)--(4.5,2)--(0,2)--(0,.5)--(3,.5)--(3,2);
 \draw (1.5,2)--(1.5,0.5);
 \draw (0,1.5)--(4.5,1.5);
 \draw (3,0.5)--(4.5,0.5)--(4.5,1);

 \draw (0.75,0.75) node {$x_3-y_1$};
 \draw (0.75,1.25) node {$x_2-y_1$};
 \draw (0.75,1.75) node {$x_1-y_1$};
 \draw (2.25,0.75) node {$x_3-y_2$};
 \draw (2.25,1.25) node {$x_2-y_2$};
 \draw (2.25,1.75) node {$x_1-y_2$};
 \draw (3.75,0.75) node {$0$}; 
 \draw (3.75,1.25) node {$x_2-x_3$};
 \draw (3.75,1.75) node {$x_1-x_3$};

\end{tikzpicture}

    \end{subfigure}
    \hfill
    \vspace{1cm}

    \begin{subfigure}[h]{0.4\textwidth}
        \centering
        \begin{tikzpicture}

 \draw (0,1)--(4.5,1)--(4.5,2)--(0,2)--(0,.5)--(3,.5)--(3,2);
 \draw (1.5,2)--(1.5,0.5);
 \draw (0,1.5)--(4.5,1.5);
 \draw (3,0.5)--(4.5,0.5)--(4.5,1);

 \draw (0.75,0.75) node {$x_3-y_1$};
 \draw (0.75,1.25) node {$x_2-y_1$};
 \draw (0.75,1.75) node {$x_1-y_1$};
 \draw (2.25,0.75) node {$x_3-y_2$};
 \draw (2.25,1.25) node {$x_2-y_2$};
 \draw (2.25,1.75) node {$x_1-y_2$};
 \draw (3.75,0.75) node {$0$}; 
 \draw (3.75,1.25) node {$x_2-x_3$};
 \draw (3.75,1.75) node {$x_1-x_3$};
 \draw (0.75,0.2) node {};

\end{tikzpicture}
 \end{subfigure}
    \hfill
    \begin{subfigure}[h]{0.1\textwidth}
        \centering   
        $\displaystyle{\lmp{L}}$
        
            \end{subfigure}
    \hfill
    \begin{subfigure}[h]{0.4\textwidth}
        \centering    
        
        \begin{tikzpicture}
 
\vspace{0.cm}

 \draw (0,1)--(4.5,1)--(4.5,2)--(0,2)--(0,.5)--(3,.5)--(3,2);
 \draw (1.5,2)--(1.5,0.5);
 \draw (0,1.5)--(4.5,1.5);

 \draw (0.75,0.75) node {$e_1-e_3$};
 \draw (0.75,1.25) node {$e_1-e_4$};
 \draw (0.75,1.75) node {$e_1-e_5$};
 \draw (2.25,0.75) node {$e_2-e_3$};
 \draw (2.25,1.25) node {$e_2-e_4$};
 \draw (2.25,1.75) node {$e_2-e_5$};
 \draw (3.75,0.75) node {$\alpha_3$}; 
 \draw (3.75,1.25) node {$e_3-e_4$};
 \draw (3.75,1.75) node {$e_3-e_5$};
 
 \draw (0.75,0.3) node {$\alpha_1$};
 \draw (2.25,0.3) node {$\alpha_2$};
 \draw (4.8,1.25) node {$\alpha_4$};
 \draw (4.8,1.75) node {$\alpha_5$};

\end{tikzpicture}

    \end{subfigure}
    \hfill    
    
        \caption{Degeneration for $\pi=[1243]$.}\label{fig:maps}
\end{figure}
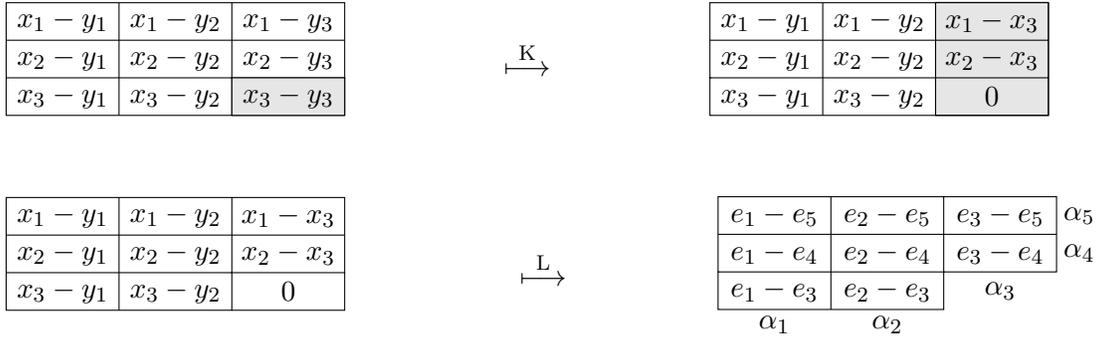

\proof We will give a linear map from $\mpy\rightarrow \P(T(\pi))$ as a composition of two maps $K$ and $L$. Let $L$ be the map given by 
\begin{align*} L(x_i)&=-e_j \text{, where } \alpha_j \text{ is the label of step } N \text{ on row } i, and \\
L(y_i)& = \begin{cases} 0 &\mbox{if } (a,i)\in Ess(\pi) \mbox{ for some }a, \\ 
-e_j & \mbox{where } \alpha_j \mbox{ is the label of step } E \mbox{ on column } i. \end{cases} 
\end{align*}
Let $K$ be the map given by 
\begin{align*}K(y_j)&=\begin{cases} x_i &\mbox{if } (i,j)\in Ess(\pi), \\ y_j & \mbox{if there is no } a \mbox{ such that } (a,j)\in Ess(\pi).
\end{cases}\\
K(x_i)&=x_i.\end{align*}
See Figure \ref{fig:maps} for an example of these maps. Now notice that if $(x_i-y_j)\in \Phi((\mathbb{P}Y_\pi)^{T^{2n-1}})$, i.e. it is the image of a fixed point, then we have that $L\circ K((x_i-y_j))=e_{k_j}-e_{k_i}$, where $\alpha_{k_i}$ is the label of step $N$ in $A$ on row $i$ and $\alpha_{k_j}$ is the label of step $E$ in $A$ on column $j$. Notice that $K$ corresponds to subtracting the vector $x_a-y_b$ from each vector $x_i-y_b$ in column $b$ whenever $(a,b)\in Ess(\pi)$. Notice also that $L$ is the map induced by the relabeling of the entries in $NW(\pi)-Ess(\pi)$ into the tree labeling. In the language of $\phi$, we then have that $(x_i-y_j)\mapsto e_{\phi^{-1}(j)}-e_{\phi^{-1}(i)}$. The root polytope of $T(\pi)$ is 
\[\textrm{ConvHull}(0,e_{\phi^{-1}(j)}-e_{\phi^{-1}(i)}\mid (i,j)\in NW(\pi)-Ess(\pi))\]
so the theorem follows.
\qed

By Propositions \ref{prop:cress} and \ref{prop:esslab} we have that for $\pi=1\pi'$ with $\pi'$ dominant the image $L\circ K(\mpy)$ is precisely the root polytope $\P(T(\pi))$.

The following Proposition allows us to describe the degeneration in Theorem \ref{thm:tortoroot} in words.

\bp The vertices of $\mpy$ corresponding to the entries in $Ess(\pi)$ form a face of $\mpy$.
\ep

\proof Let $Ess(\pi)=\{(a_1,b_1),\ldots,(a_k,b_k)\}$ and consider the linear functional $\ell$ defined by
\begin{equation*}
\ell\left(\sum_i s_ix_i+\sum_jt_jy_j\right) :=\sum_{i=1}^k is_{a_i} +\sum_{j=1}^k(j-k-1)t_{b_j}.
\end{equation*}
Note that since $(a_i,b_j)\notin L(\pi)$ if $i<j$, then the maximum value $\ell$ can attain on $\mpy$ is $k+1$. Moreover, this value is only attained for the vectors $x_a-y_b$ with $(a,b)\in Ess(\pi)$.
\qed

The degeneration $L\circ K:\mpy\rightarrow \P(T(\pi))$ consists of contracting the face of $\mpy$ corresponding to $Ess(\pi)$ to a point and moving this point to the origin while tweaking the vertices of $\mpy$ that are of the form $\frac{1}{2}(x_i-y_j)$ where $(i,j)$ is north of an entry of $Ess(\pi)$.

\subsection{Triangulating $\mpy$ and geometric realization of subword complexes} Next we show that the image of the  canonical triangulation of $\P(T(\pi))$ for $\pi=1\pi'$, with $\pi'$ dominant, under the linear map $L\circ K$ is a triangulation of $\Phi(\mathbb{P}(Y_\pi))$, which is yet another way to geometrically realize the pipe dream complex $PD(\pi)$ for these permutations. 

\bt \label{can} Let $\Delta_1, \ldots, \Delta_k$ be the top dimensional simplices in the  canonical triangulation of $\P(T(\pi))$  for $\pi=1\pi'$, where $\pi'$ is dominant. Then $P_i:=(L\circ K)^{-1}(\Delta_i)$, $i \in [k]$, are the top dimensional simplices in  a triangulation of $\Phi(\mathbb{P}(Y_\pi))$ which we call its \textbf{canonical triangulation}.
\et

\proof Since $L\circ K(\mpy)=\P(T(\pi)),$ it follows that $\cup_{i \in [k]} P_i=\Phi(\mathbb{P}(Y_\pi))$. We will show that the $P_i$ are interior disjoint simplices that induce a triangulation of $\Phi(\mathbb{P}(Y_\pi))$. 
 
Note that each vertex of $\P(T(\pi))$ except $0$ has  a  unique vertex of $\Phi(\mathbb{P}(Y_\pi))$ that maps to it under $L\circ K$ and there are $|Ess(\pi)|$ many vertices of $\Phi(\mathbb{P}(Y_\pi))$ mapping to $0$. 
For $v\neq 0$ a vertex of  $\P(T(\pi))$, let $w(v)$ be the vertex of $\Phi(\mathbb{P}(Y_\pi))$ that maps to it under $L\circ K$. Let $q_1, \ldots, q_{|Ess(\pi)|}$ be the vertices of $\Phi(\mathbb{P}(Y_\pi))$ mapping to  $0$ under $L\circ K$. We have that  \be \label{pi} P_i=ConvHull(q_1, \ldots, q_{|Ess(\pi)|}, w(v) \mid  v \text{ is a nonzero vertex of } \Delta_i),\ee 
since $(L\circ K)^{-1}(0)=ConvHull(q_1, \ldots, q_{|Ess(\pi)|})$ and $(L\circ K)^{-1}(v)=w(v)$ for $ v$ a nonzero vertex of $\P(T(\pi))$. 
Notice that $\dim(\Delta_i)=\dim(\P(T(\pi)))$ which equals to the number of edges of $T(\pi)$ which is $r+c-1-|Ess(\pi)|$ and so it follows that $P_i$ is the convex hull of $r+c-1$ points.
It remains to show that the dimension of each $P_i$ equals the dimension of $\mpy$.
This can be done by fixing a vertex $u$ of $\Delta_i$ and showing that the set $\{w(u)-q_1,\ldots,w(u)-q_{|Ess(\pi)|},w(v)-w(u) \mid  v \text{ is a nonzero vertex of } \Delta_i\}$ is linearly independent. We leave the details of this argument to the interested reader.

\qed

\medskip

Denote by $\C(\pi)$ the core of the pipe dream complex $PD(\pi)$. Let $\C^i(\pi)$ be the core $\C(\pi)$ coned over $i$ times. 

\bt \label{spec} The canonical triangulation of $\Phi(\mathbb{P}(Y_\pi))$, for $\pi=1\pi'$, with $\pi'$ dominant, is a geometric realization of $\C^{|Ess(\pi)|+1}(\pi).$
\et

\proof By \cite[Theorem 1]{us} the canonical triangulation of $\P(T(\pi))$ is a geometric realization of $\C^2(\pi)$.  For $\pi$ as in the theorem, we have that $dom(\pi)$ is empty and $L(\pi)$ is a partition with $r$ rows and $c$ columns. Then by equation \eqref{p} and Lemma  \ref{cor-dim} we get that $dim \Phi(\mathbb{P}(Y_\pi))=r+c-2$.   On the other hand, the dimension of $\P(T(\pi))$ equals the number of edges of $T(\pi)$, which is $r+c-1-|Ess(\pi)|$. Thus, $dim(\P)-dim (\P(T(\pi)))=|Ess(\pi)|-1$. Together with Theorem \ref{can} we obtain the theorem. \qed

\section*{Acknowledgements} We are grateful to Allen Knutson for many valuable discussions. We thank Maksim Maydanskiy for his interest and indepth comments on our research. We also thank Vic Reiner for several interesting conversations. 

\medskip

\bibliography{biblio-kir}
\bibliographystyle{alpha}

\end{document}